\documentclass[11pt,
]{article}

\usepackage[arrow,matrix,curve]{xy}

\usepackage{amssymb, latexsym, amsmath, amscd, array, amsthm}
\usepackage{mathrsfs}
\usepackage{mparhack}

\input{definable.sty}

\begin{document}

\title 
{Minimal axiomatic frameworks for definable hyperreals with transfer}

\author {Frederik S.\ Herzberg \thanks{ Center for Mathematical
Economics (IMW) and Institute for Interdisciplinary Studies of Science
(I$^2$SoS), at Bielefeld University, Universit\"atsstra\ss{}e 25,
D-33615 Bielefeld, Germany, and Munich Center for Mathematical
Philosophy (MCMP), Ludwig Maximilian University of Munich,
Geschwister-Scholl-Platz 1, D-80539 Munich, Germany.  {\tt
fherzberg@uni-bielefeld.de}.  }  \and Vladimir Kanovei \thanks{ IITP,
Moscow, and MIIT, Moscow, Russia.  Partial support of ESI at Vienna,
during the visit in December 2016, acknowledged.  Partial support of
grant RFBR 17-01-00705 acknowledged.  {\tt kanovei@googlemail.com}}
\and Mikhail Katz \thanks{Bar Ilan University, Ramat Gan 5290002,
Israel.  Supported in part by Israel Science Foundation grant 1517/12.
{\tt katzmik@macs.biu.ac.il}} \and Vassily Lyubetsky \thanks{IITP,
Moscow, Russia. Partial support of grant RSF 14-50-00150 acknowledged.
{\tt lyubetsk@iitp.ru}} }

\date{\today}
\maketitle

\begin{abstract}
We modify the definable ultrapower construction of Kanovei and Shelah
(2004) to develop a~$\ZF$-definable extension of the continuum with
transfer provable using countable choice only, with an additional mild
hypothesis on well-ordering implying properness.  
Under the same assumptions, we
also prove the existence of a definable, proper elementary extension
of the standard superstructure over the reals.

\end{abstract}

\subsection{Introduction}
\label{sint}

The usual ultrapower construction of a hyperreal field~$\dR^\om/U$ is
not functorial (in the category of models of set theory) due to its
dependence on a choice of a free ultrafilter~$U$, which can be
obtained in~$\ZFC$ only as an application of the axiom of choice, but
not as an explicitly definable set-theoretic object.  Kanovei and
Shelah \cite{KS} developed a functorial alternative to this, by
providing a construction of a definable hyperreal field, which we
refer to below as \emph{the KS construction}.

The KS construction was analyzed in \cite{Ke76}, Section 1G of the
Online 2007 edition, and generalized in many ways in
\cite{herz,klyub,kusp}, \cite[Chapter 4]{krbook} among others.  We
give \cite{keis} as the source of the problem of a uniquely definable
nonstandard real line, and \cite{lux}, \cite{SL}, \cite{Da77},
\cite{HL} as basic references in nonstandard matters.

Nonstandard analysis is viewed by some as inherently non-constructive.
One of the reasons is that nonstandard models are typically presented
in terms of an unspecified choice of a free ultrafilter, which makes
the resulting ultrapower hopelessly non-definable.  In fact, as
Luxemburg \cite{lux} observed, if there is a non-standard model of the
reals, then there is a free ultrafilter on the natural numbers~$\om$.
The observation that elements of a nonstandard extension
${}^{\ast}\!\!A$ correspond to ultralters on~$A$ was first exploited
in detail by Luxemburg \cite{Lu69}.

To circumvent the unspecified choice of a free ultrafilter, the KS
construction starts with the collection of ultrafilters~$U$ on~$\om$
%(both free and principal) 
parametrized by surjective maps from a suitable ordinal onto such
ultrafilters~$U$.  Such maps are ordered lexicographically.  This
generates a definable linear ordering of ultrafilters in which each of
them is included in many copies.  The tensor product is applied to
merge the ultrafilters into a definable ultrafilter in the algebra of
finite support product sets.

Thus, the KS construction can be viewed as a functor which, 
given a model of set theory, produces a definable extension 
of the reals in the model.

The KS construction in \cite{KS}
(as well as its modifications as in \cite{herz,kusp}) 
was originally designed to work in Zermelo--Fraenkel 
set theory~$\ZFC$ with choice.  
However for it {\sl prima facie\/} to yield the  
expected result, it is sufficient to assume the 
wellorderability of the reals. 
%the continuum~$\cP(\om)=\ans{X:X\sq\om}$ is wellorderable.
%  
Let~$\wor$ be the following statement: the
continuum~$\dn=\ens{X}{X\sq\om}$ is wellorderable.  Thus it emerges
that the theory~$\ZF + \wor$ is sufficient for the KS construction to
yield a definable proper elementary extension of the reals.
 
The goal of this note is to weaken this assumption.

\subsection{The result}
\label{res}

Consider the following two consequences of the axiom of 
choice in~$\ZF$: 
\bde
\item[\cac:] countable \AC{} for sets of reals, 
that is, any sequence~$\sis{X_n}{n<\om}$ 
of sets~$\pu\ne X_n\sq \dR$ admits a
choice function;

\item[\wob:] there exists a free ultrafilter over~$\om$ with a
wellorderable base.  (A set~$B\sq U$ is a {\it base} of an ultrafilter
$U$ over~$\om$, if and only if there is no ultrafilter~$U'\ne U$
over~$\om$ with~$B\sq U'$.  In such case we write~$U=\uob B$.)
%
%or in short, a free {\it wob-ultrafilter}.  
\ede

\bte[$\ZF$]
\label{t}
There exists an extension\/~$\adR$ of the reals\/~$\dR$, such that
both\/~$\adR$ and a canonical embedding\/~$x\longmapsto \ax$ from\/
$\dR$ into\/~$\adR$ are presented by explicitly definable
set-theoretic constructions, and in addition$:$ \ben \renu \itla{t2}
$\cac$ implies that\/~$\adR$ is an elementary extension, in the sense
of the language\/~$\cL(\dR)$ with symbols for all finitary relations
on\/~$\dR\;;$

\itla{t3}
$\wob$ implies that\/~$\adR$ is a proper extension of\/ 
$\dR$, 
%therefore 
containing infinitesimals and infinitely 
large numbers.
\een
\ete

It follows by \ref{t2} that, instead of~$\wor$, the axiom~$\cac$ can
be used to establish elementarity.  It emerges that proving the
transfer principle for the definable extension requires no more choice
than proving, for instance, the~$\sigma$-additivity of the Lebesgue
measure; see \cite{KK}.  Similarly, by \ref{t3},~$\wob$ successfully
replaces~$\wor$ in the proof of properness.

Quite obviously~$\wor$ implies~$\cac$ and~$\wob$ in~$\ZF$.  The
failure of the inverse implication is dealt with in \ref{s2}.1 below.

The proof of Theorem~\ref{t} appears in Section~\ref{s1}.  We also
show, in \ref{s2}.3, how the theorem can be generalized in order to
obtain even a nonstandard superstructure over~$\adR$.

\subsection{What it takes: array of ultrafilters}
\label{s1}

Let an {\it array of ultrafilters\/} be any sequence 
$\sis{D_a}{a\in A}$, where~$A=\ang{A,<_A}$ is a 
linearly ordered set and each~$D_a$ is an ultrafilter 
over~$\om$. 

\bpro
[in~$\ZF+\cac$]
\label{kst}
Assume that\/~$\sis{D_a}{a\in A}$ is a 
definable array of ultrafilters over\/~$\om$, 
with at least one free ultrafilter\/~$D_{a_0}$. 
Then there is a definable\/ 
{\rm(as in Theorem~\ref{t})} proper extension\/~$\adR$ of\/ 
$\dR$, elementary w.\,r.\,t.\ the language\/ 
$\cL(\dR)$ containing 
%(symbols for) 
all finitary relations on\/~$\dR$. 
\epro
\bpf[sketch, based on the proof in \cite{KS}]
The following is defined:
\bit
\item[$-$]
the {\it index set\/} 
$I=\om^A=\ens{x}{x\,\text{ is a map }\,A\to\om}$;

\item[$-$] the algebra~$\cX=\cX(A)$ of {\em finite-support subsets} of
$I=\om^A$, so that a set~$X\sq\om^A$ is in~$\cX$ if and only if there
is a finite~$u\sq A$ such that
\dm
\kaz x \zi y \in \om^A\; \big( {x \res u= y \res u}
\limp{({x \in X}\eqv{ y \in X})} \big)\,;
\dm

\item[$-$] the collection~$F=F(A)$ of {\em finite-support functions}
$f:I\to\dR$, so that~$f:I\to\dR$ belongs to~$F$ if and only if there
is a finite set~$u\sq A$ such that
\dm
\kaz x \zi y \in \om^A\; \big( {x \res u= y \res u}\limp
{(f(x)=f(y))} \big)\,.
\dm
\eit
The {\it tensor, or Fubini product\/}
$D=\bigotimes_{a\in A}D_a$ consists then of all sets 
$X\sq I$ such that for a finite subset 
$u=\{a_1<_A\dots<_A a_n\}\sq A$ we have: 
$$
D_{a_n}k_n\dots D_{a_2}k_2\:D_{a_1}k_1\: 
(\ang{k_1,\dots,k_n}\inu X)\,,
$$
where  
$\ang{k_1,\dots,k_n}\inu X$ 
means that every~$x\in I$ satisfying 
$x(a_1)=k_1,\dots,\linebreak[0]x(a_n)=k_n$ 
belongs to~$X$, and 
$D_a k\,\Phi(k)$ means that the set 
$\ens{k}{\Phi(k)}$ belongs to~$D_a$.
It turns out that~$D$ is an ultrafilter in the algebra~$\cX$, which
allows to define the {\it ultrapower\/}~$\adR=F/D=\ens{\kla f}{f\in
F}$, where~$\kla f=\ens{g\in F}{f=^D g}$ and~$f\eqD g$ means that
$\ens{x \in I}{f(x )=g(x )} \in D$.  All finitary relations in
$\cL(\dR)$ extend to~$\adR$ naturally.

In addition, we send every real~$r$ to the equivalence class
$\arr=\kla{c_r}$ of the constant function~$c_r\in F$ with value~$r$.
The axiom~$\cac$ is strong enough to support the ordinary proof of the
\los\ lemma, and hence~$r\mto\arr$ is an elementary embedding in the
sense of the language~$\cL(\dR)$.
To prove that the embedding is proper, make use of 
the assumption that at least one of~$D_a$ is a free 
ultrafilter.
Finally, the extension~$\adR$ is definable since the given array of
ultralters~$\sis{D_a}{a\in A}$ is definable by hypothesis.  \epf

\bpf[Theorem~\ref{t}]  
To define a suitable array of ultrafilters, 
let~$\vt$ be the least ordinal such that for 
any wellorderable set
$Z\sq\dR$ there is a surjective map~$a\colon\vt\na Z$.
Let~$A$ consist of all maps 
$a\colon\vt\to\cP(\om)$  such that the set
$B_a=\ran a=\ens{a(\ga)}{\ga<\vt}$ is a base of an
%nontrivial (\ie, containing only infinite sets) 
ultrafilter on~$\om$, and let~$D_a=\uob{B_a}$ be this ultrafilter.
The set~$A$ is ordered lexicographically:~$a<_A b$ if and only if
there exists an ordinal~$\ga<\vt$ such that~${a\res\ga}={b\res\ga}$
and~$a(\ga)<b(\ga)$ in the sense of the lexicographical linear order
$<$ on~$\cP(\om)$.  Then~$\sis{D_a}{a\in A}$ is a definable array of
ultrafilters.  Assuming \wob, it contains at least one free
ultrafilter\/~$D_{a}$, and we apply Proposition \ref{kst}.  \epf

\subsection{Remarks}
\label{s2}

Here we add some related remarks, starting with a model
of~$\ZF+\acc+\wob{}$ in which the continuum is not wellorderable.
This demonstrates that Theorem~\ref{t} is an actual strengthening of
the key result of \cite{KS}.\vom

{\ubf\arabic{subsection}.1.  Separating~$\wob+\acc$ from~$\wor$.}
\label{s41}
Pincus and Solovay conjectured in \cite[p.\;89]{ps} that iterated
Sacks extensions may be useful in the construction of choiceless
models with free ultrafilters.  Working in this direction, we let
$\gM$ be an~$\omega_1$-iterated Sacks extension of~$\bL$, the
constructible universe, as in~\cite{bl}.
Let~$\gN$ be the class of all sets hereditarily   
definable from an~$\om$-sequence of ordinals in~$\gM$.  
Then \wob\ is true in~$\gN$ since some free ultrafilters 
in~$\bL$ 
(basically, all selective ultrafilters) 
remain ultrafilter bases in~$\gM$ and in~$\gN$ by
\cite[Section\;4]{bl}, and \acc\ is true as well 
(even the principle of dependent choices \dc\ holds). 
Meanwhile,~$\wor$ fails in~$\gN$ 
($\dn$ is not wellorderable) 
by virtue of arguments, based on the homogeneous 
structure of the Sacks forcing, and similar to those 
used in the classical studies of the 
choiceless Solovay model~$S'$ as in 
\cite[Part III, proof of Theorem\;1]{sol}.

If we let~$\gM$ be an~$\omega_2$-iterated Sacks extension of~$\bL$,
then the class~$\gN$ of all sets that are hereditarily definable from
an~$\omi$-sequence of ordinals in~$\gM$, still will be a model of
$\wob+\neg\wor$, in which even~$\dc_{\omi}$ and~$\AC_{\omi}$ hold
instead of the simple \dc\ and \acc.  Longer iterations make little
sense here as each further Sacks real collapses all smaller cardinals
down to~$\omi$.

We know nothing about any model of~$\wob+\cac$ in 
which~$\wor$ fails, 
different from the ones just described. 
(However see \ref{s2}.5 below.) 
This can be a difficult problem, yet not uncommon in 
studies of choiceless models.\vom

{\ubf\ref{s2}.2. 
Keisler-style representation.} 
Keisler's influential monograph \cite{Ke76} contains 
(in Section 1G) a somewhat modified exposition of 
the construction of a definable
nonstandard extension of \cite{KS}, by an explicit 
{\it amalgamation\/} of all ultrapowers of~$\dR$ 
via different ultrafilters on~$\om$ into one large 
hyperreal field.  
A similar Keisler-style modification of the construction 
readily works in the~$\ZF+\cac+\wob$ setting.
\vom

{\ubf\ref{s2}.3. 
Superstructure over~$\adR$.}
Let~$V(\dR)=\bigcup_{n\geq0}V_n(\dR)$ be the superstructure over the
reals, where~$V_0(\dR)=\dR$
and~$V_{n+1}(\dR)=V_n(\dR)\cup\cP(V_n(\dR))$ for all~$n$, see
\cite[Section\;4.4]{ck}.  To build a nonstandard superstructure
over~$\adR$ as in Section~\ref{s1}, we let~$F_n$ be the set of all
functions~$f:\ia\to V_n(\dR)$ of finite support, and then define the
ultrapower~$\sV_n(\dR)=F_n/D$ and the elementary embedding~$x\mto\ax$
from~$V_n(\dR)$ to~$\sV_n(\dR)$ as above.  (And we need~$\acc$ for
subsets of~$V_n(\dR)$ to prove the elementarity.)  Then each element
of~$\sV_n(\dR)$ can be identified with a certain subset of
$\sV_{n-1}(\dR)$ or an element of~$\sV_{n-1}(\dR)$, so that each
$\sV_n(\dR)$ emerges as a subset of~$V_{n}(\adR)$.  This completes the
nonstandard superstructure construction under~$\wob+\acc$.\vom

{\ubf\ref{s2}.4.  
Another definable choiceless ultrapower.}
Consider {\em the basic Cohen model\/}~$\bL(A)$, 
obtained by adding a set
$A=\ens{a_n}{n<\om}$ of Cohen generic reals~$a_n$ to~$\bL$,  
\cite[5.3]{jac}. 
(Not to be confused with {\em the Feferman model\/} 
\cite%[4.12]
{fef}, 
adding all~$a_n$ but not~$A$.) 
%and does not have free ultrafilters over~$\om$,)  
The set~$A$ belongs to~$\bL(A)$ but the map~$n\mto a_n$ 
does not. 
$\acc$ badly fails in~$\bL(A)$ as~$A$ is an infinite 
Dedekind finite set.
Yet 
%it is a more peculiar property of this model that it 
$\bL(A)$ contains a free ultrafilter~$U$ over~$\om$,  
see \cite{rep} for a short proof.

Let~$\adR=\dR^\om/U$ be the associated ultrapower.  Then~$\adR$ {\ubf
is not} an elementary extension of~$\dR$ in the full relational
language~$\cL(\dR)$ as in Theorem~\ref{t}, since the formula ``$\kaz
n\in\om\:\,\sus x$($x$ codes an~$n$-tuple of elements of~$A$)'' is true
for~$\dR$ but false for~$\adR$.
%
%(See Howard~\cite{how} for a general argument.)
%On the other hand, a simple argument shows that 
%
However~$\adR$ {\bf is} an elementary extension of~$\dR$ with respect
to the sublanguage~$\cL'(\dR)$ of~$\cL(\dR)$, containing only
real-ordinal definable finitary relations on~$\dR$.  Note that
$\cL'(\dR)$ is a sufficiently rich language to enable an adequate
development of nonstandard real analysis.

Both~$U$ and~$\adR$ are definable in~$\bL(A)$ by a 
set theoretic formula with the only parameter~$A$.
And this is probably all we can do in~$\bL(A)$
since the model contains no real-ordinal 
definable elementary extensions of~$\dR$.
\vom

{\ubf\ref{s2}.5.  A possible \wob\ model.}  

One may want to extend~$\bL(A)$ as in \ref{s2}.4 by a~$P(U)$-generic
real~$c=c_0\in\dn$, where~$P(U)$ is the Mathias forcing with infinite
conditions in~$U$.  If~$\bL(A)[c_0]$ happens to have an
$\ans{A,a_0}$-definable ultrafilter~$U_1$ over~$\om$ with~$U\sq U_1$
then let~$c_1\in\dn$ be a~$P(U_1)$-generic real over~$\bL(A)[c_0]$.
Extending this forcing iteration as in \cite [A10 in Chapter 8] {kun},
one may hope to get a final extension of~$\bL(A)$ with a wellordered
ultrafilter base~$\ens{c_\xi}{\xi<\omi}$ but with~$A$ still not
wellorderable.

{\ubf\ref{s2}.6.  Least cardinality.}  

What is the least possible cardinality of a definable hyperreal field?
A rough estimate for the general definable extension in \cite{KS}
under AC yields $\le {\exp^3(\aleph_0)}$.  As for the definable
extension $\adR$ in Section~\ref{s1} of this paper, if the ground set
universe is the $\omega_2$-iterated Sacks extension of $\bL$ as
in~\ref{s41}.1 then $\text{card}(\adR)=2^{\aleph_0}=\aleph_2$, which
is minimal.

\subsection{Conclusions}
\label{conc}

Analysis with infinitesimals presupposes the existence of an extended
mathematical universe which, in the tradition of Robinson and Zakon
\cite{RZ}, is typically understood as an extended superstructure over
the reals, although for some basic applications an extension of the
set of reals suffices.  Even for certain more sophisticated
applications, it is enough for this extension of the mathematical
universe to satisfy the Transfer Principle, which means that it is an
elementary extension in the sense of model theory.

We have shown that one can find definable extensions of both the set
of reals and the superstructure over the reals; more precisely, our
extensions are definable by purely set-theoretic means without
recourse to well-ordering, and have the following properties: (I) one
can prove the \emph{Transfer Principle} for such extensions from
Zermelo--Fraenkel set theory with merely Countable Choice; (II) the
existence of infinitesimals and infinitely large numbers in those
extensions follows from a mild well-ordering assumption.
% (viz., the existence of free
%ultrafilters on~$\omega$ with a wellorderable base, \wob).  

The property of \emph{countable saturation},
% of those models, 
important for some advanced applications, 
% of nonstandard methods, 
is not asserted, but can be achieved by the 
$\omega_1$-iteration of the given
extension construction, as described in \cite[Section 4]{KS}.

Our results may be of interest to practitioners working with fragments
of nonstandard analysis.  For instance, the Transfer Principle plus
the existence of an infinitely large integer is all that is required
to develop Edward Nelson's \cite[p.~30]{N07} \emph{minimal nonstandard
analysis} or the related \emph{minimal Internal Set Theory}
\cite[pp.\;3, 4, 104]{herz+}.  Such fragments of nonstandard analysis
have the potential for application in diverse fields, ranging from
stochastic calculus and mathematical finance to theoretical quantum
mechanics \cite{herz+}.

{\ubf Acknowledgements.}  

The authors are grateful to the anonymous JSL referee whose remarks
and suggestions helped improve the article.

\renek{\refname}
\end{document}